# Looking for good Hofmeister and Braunschadel bases

## 0 History

| | |
|---|---|
| 28-10-92 | Document started |
| 01-11-92 | First version completed |
| 24-01-93 | New sections 6 and 7 added to record further results (for t=21 to 58), and issues arising from them. |
| | Summary table of the "best" bases added as Appendix 1. |
| 14-03-93 | Note on extended search (which confirms previous results) added. |
| 29-09-14 | Abstract and Reference sections added, and internal references updated, to prepare the document for ArxIv publication. Some further explanatory material added, as well as a new Appendix containing a transcript of Selmer's "private communication" to me dated 1992/1993) (see Selmer, E. S., [8]). |
| | Other than this additional explanatory material there are no changes to the content of the paper as originally prepared in March 1993. |

## Abstract

$A_k$ = (1, $a_2$, ... $a_k$} is an h-basis for n if every positive integer not exceeding n can be expressed as the sum of no more than h values $a_i$. An "extremal" h-basis $A_k$ is one for which n is as large as possible. Computing extremal bases has become known as the "global" Postage Stamp Problem.

This paper describes the author's early attempts to identify extremal Hofmeister and Braunschadel bases for large t (where t = 12h + r for some 0 <= r <= 11), and also includes a transcription of Prof. Selmer's correspondence with the author about this work which has not been published before. See (Selmer, E.S., Private Communication, [8], 1992/1993).

## 1 Early work

In August 1989 a colleague and personal friend Hector Prior suggested formulae of the following kind for maximal d=4 sets:

$$\begin{aligned}
h &= 12t + r \\
a_4 &= (k_{43}t + c_{43})a_3 + (k_{42}t + c_{42})a_2 + (k_{41}t + c_{41}) \\
a_3 &= (k_{32}t + c_{32})a_2 + (k_{31}t + c_{31}) \\
a_2 &= (k_{21}t + c_{21})
\end{aligned}$$

with the $k_{ij}$ independent of both t and r, and the $c_{ij}$ dependent on r only. Using the computed maximal sets for h <= 76, he suggested the following values for the $k_{ij}$:

$$(k_{21}, k_{31}, k_{32}, k_{41}, k_{42}, k_{43}) = (9, 4, 3, 6, 2, 2)$$

I then investigated formulae of this kind by computer, and it soon became apparent that the value for $k_{41}$ should be 7, rather than 6. I now concentrated on identifying the $c_{ij}$, and developed a program (called exp21c) which took as input:

- a fixed value of r

- an initial value of t
- a range of values for each coefficient $c_{ij}$

and then, for successive values of t, considered each set in turn, printing out those whose cover was at least as good as that of some "good" set (as defined by a particular set of coefficients $c_{ij}$). This program was, in fact, the first of my programs to employ the "difficult target" idea (see Challis, M.F., [3]) and I was able to investigate ranges of $c_{ij}$ right up to t=20.

At this time, I believed that there must exist a "finite set of formulae" to describe the d=4 maximal sets that would be valid for all sufficiently large h, and also that we were close to identifying those formulae; clearly, the cycle length would be 12, and almost certainly we now had the correct values for the coefficients $k_{ij}$; the values I proposed for the $c_{ij}$ in September 1989 were:

Table 1

| r | $c_{21}$ | $c_{31}$ | $c_{32}$ | $c_{41}$ | $c_{42}$ | $c_{43}$ |
|---|---|---|---|---|---|---|
| 0 | 1 | 0 | 0 | 0 | 0 | 0 |
| 1 | 1 | 0 | 2 | 1 | 1 | 0 |
| 2 | 1 | 0 | 2 | 1 | 1 | 0 |
| 3 | 3 | 1 | 2 | 2 | 1 | 1 |
| 4 | 3 | 1 | 2 | 2 | 1 | 1 |
| 5 | 3 | 1 | 2 | 2 | 1 | 1 |
| 6 | 3 | 1 | 2 | 2 | 1 | 1 |
| 7 | 8 | 4 | 1 | 7 | 1 | 0 |
| 8 | 8 | 4 | 1 | 7 | 1 | 0 |
| 9 | 8 | 4 | 1 | 7 | 1 | 0 |
| 10 | 11 | 6 | 1 | 10 | 1 | 0 |
| 11 | 11 | 5 | 2 | 9 | 2 | 0 |

**2 Braunschadel and Hofmeister bases**

This is where the matter remained until late 1991 when I sent a copy of my results to Professor Ernst Selmer, recently retired from the University of Bergen. By this time, I had applied the "difficult target" technique to the d=4 case (in the program gen27c), and had calculated extremal sets as far as h=122. Selmer reviewed these results, and, in a letter dated January 2nd 1992, noted that two of them - for h=107 and 119 (ie r=11 for t=8, 9) - did not fit into the formula above; instead, he showed that they conformed to a very similar formula in which just one coefficient $k_{ij}$ was different: $k_{31} = 2$ instead of 4.

Selmer had recently been working on the concept of "dual" or "associate" bases (Selmer, E. S., [9]) and this new basis is the dual of the original one; it had first been discovered and investigated by Braunschadel (Braunschadel, R., [1]). It also turned out that formulae of the original kind had been discovered by Hofmeister and had been thoroughly investigated by Mossige in his thesis submiited in 1986 (Mossige, Svein, [6]), leading to the discovery of bases with an asymptotic coefficient of 2.008... - a most surprising result.

For these reasons, I now refer to bases for h = 12t + r as follows:

$(k_{21}, k_{31}, k_{32}, k_{41}, k_{42}, k_{43}) = (9, 4, 3, 7, 2, 2)$ *Hofmeister* bases

$(k_{21}, k_{31}, k_{32}, k_{41}, k_{42}, k_{43}) = (9, 2, 3, 7, 2, 2)$ *Braunschadel* bases

Further results for d=4 (up to h=158 at the time of writing) suggest that all maximal sets are either of the Hofmeister or of the Braunschadel form, and so I returned to exp21 to see what coefficients $c_{ij}$ were best for larger values of t; in particular, there seemed to be some evidence that Braunschadel bases were beginning to take over from Hofmeister for large t.

The results are summarised in Tables 2 and 3 below, details of which I sent to Selmer on April 5th 1992. Although I shall talk of the *best* Hofmeister or Braunschadel basis for a given value of t, it is important to remember that I have not completed an exhaustive search: the search is only over those sets defined by coefficients $c_{ij}$ that lie within the ranges given in the tables.

Table 2
Probably the best Hofmeister bases for 5<=t<=20

| | | best coefficients | | | | | | range of coefficients checked | | | | | |
|---|---|---|---|---|---|---|---|---|---|---|---|---|---|
| r | $t_{range}$ | $c_{21}$ | $c_{31}$ | $c_{32}$ | $c_{41}$ | $c_{42}$ | $c_{43}$ | $c_{21}$ | $c_{31}$ | $c_{32}$ | $c_{41}$ | $c_{42}$ | $c_{43}$ |
| 0 | [4-5] | 2 | 1 | 0 | 1 | 0 | 1 | | | | | | |
| | [6-20] | 1 | 0 | 0 | 0 | 0 | 0 | [-9,11] | [-5,5] | [-3,3] | [-10,10] | [-3,3] | [-2,2] |
| 1 | [5-20] | 1 | 0 | 2 | 1 | 1 | 0 | [-9,11] | [-5,5] | [-1,5] | [-9,11] | [-2,4] | [-2,2] |
| 2 | [5-6] | 2 | 1 | 1 | 1 | 1 | 1 | | | | | | |
| | [7-20] | 1 | 0 | 2 | 1 | 1 | 0 | [-9,11] | [-5,5] | [-1,5] | [-9,11] | [-2,4] | [-2,2] |
| 3 | [1-20] | 3 | 1 | 2 | 2 | 1 | 1 | [-7,13] | [-4,6] | [-1,5] | [-8,12] | [-2,4] | [-1,3] |
| 4 | [2-20] | 3 | 1 | 2 | 2 | 1 | 1 | [-7,13] | [-4,6] | [-1,5] | [-8,12] | [-2,4] | [-1,3] |
| 5 | [4-20] | 3 | 1 | 2 | 2 | 1 | 1 | [-7,13] | [-4,6] | [-1,5] | [-8,12] | [-2,4] | [-1,3] |
| 6 | [5-20] | 3 | 1 | 2 | 2 | 1 | 1 | [-7,13] | [-4,6] | [-1,5] | [-8,12] | [-2,4] | [-1,3] |
| 7 | [2-11] | 7 | 3 | 2 | 5 | 1 | 2 | | | | | | |
| | [12-20] | 8 | 4 | 1 | 7 | 1 | 0 | [-2,18] | [-1,9] | [-2,4] | [-3,17] | [-2,4] | [-2,2] |
| 8 | [1-16] | 7 | 3 | 3 | 5 | 2 | 2 | [-3,17] | [-2,8] | [0,6] | [-5,15] | [-1,5] | [0,4] |
| | [17-20] | 8 | 4 | 1 | 7 | 1 | 0 | [-3,17] | [-2,8] | [0,6] | [-5,15] | [-1,5] | [0,4] |
| 9 | [1-20] | 7 | 3 | 3 | 5 | 2 | 2 | [-3,17] | [-2,8] | [0,6] | [-5,15] | [-1,5] | [0,4] |
| 10 | [4-19] | 7 | 3 | 3 | 5 | 2 | 2 | [-3,17] | [-2,8] | [0,6] | [-5,15] | [-1,5] | [0,4] |
| | [20] | 11 | 6 | 1 | 10 | 1 | 0 | [-3,17] | [-2,8] | [0,6] | [-5,15] | [-1,5] | [0,4] |
| 11 | [2-17] | 10 | 4 | 3 | 7 | 2 | 2 | [0,20] | [-1,9] | [0,6] | [-3,17] | [-1,5] | [0,4] |
| | [18-20] | 11 | 5 | 2 | 9 | 2 | 0 | [0,20] | [-1,9] | [0,6] | [-3,17] | [-1,5] | [0,4] |

Notes: This table is valid for all t>=5 (h>=60), and, in some cases, for smaller values of t; such cases are noted in the $t_{range}$ column.
Where no ranges are defined, this is because the given set is known to be the *extremal* basis for those values of t.

Table 3
Probably the best Braunschadel bases for 12<=t<=20

| | | best coefficients | | | | | | range(s) of coefficients checked | | | | | |
|---|---|---|---|---|---|---|---|---|---|---|---|---|---|
| r | $t_{range}$ | $c_{21}$ | $c_{31}$ | $c_{32}$ | $c_{41}$ | $c_{42}$ | $c_{43}$ | $c_{21}$ | $c_{31}$ | $c_{32}$ | $c_{41}$ | $c_{42}$ | $c_{43}$ |
| * 0 | [12-20] | 2 | 2 | -1 | 3 | -1 | 0 | [-5,15] | [-10,10] | [-5,5] | [-5,15] | [-5,5] | [-2,2] |
| 1 | [12-15] | 0 | -1 | 2 | -1 | 0 | 1 | [-5,15] | [-10,10] | [-5,5] | [-5,15] | [-5,5] | [-2,2] |
| | [16-20] | 2 | 2 | -1 | 3 | -1 | 0 | [-5,15] | [-10,10] | [-5,5] | [-5,15] | [-5,5] | [-2,2] |
| 2 | [12] | 0 | -1 | 2 | -1 | 0 | 1 | [-5,15] | [-10,10] | [-5,5] | [-5,15] | [-5,5] | [-2,2] |
| | [13-20] | 5 | 3 | -1 | 6 | -1 | 0 | [-5,15] | [-10,10] | [-5,5] | [-5,15] | [-5,5] | [-2,2] |
| 3 | [12-20] | 5 | 3 | 0 | 6 | 0 | 0 | [-5,15] | [-10,10] | [-5,5] | [-5,15] | [-5,5] | [-2,2] |
| 4 | [12-17] | 3 | 0 | 3 | 2 | 1 | 1 | [-5,15] | [-10,10] | [-5,5] | [-5,15] | [-5,5] | [-2,2] |
| | [18-20] | 5 | 3 | 0 | 6 | 0 | 0 | [-5,15] | [-10,10] | [-5,5] | [-5,15] | [-5,5] | [-2,2] |
| 5 | [12-20] | 3 | 0 | 3 | 2 | 1 | 1 | [-5,15] | [-10,10] | [-2,8] | [-5,15] | [-5,5] | [-2,2] |
| 6 | [12-20] | 3 | 0 | 3 | 2 | 1 | 1 | [-5,15] | [-10,10] | [-2,8] | [-5,15] | [-5,5] | [-2,2] |
| 7 | [12-20] | 8 | 4 | 1 | 9 | 1 | 0 | [-5,15] | [-10,10] | [-2,8] | [-5,15] | [-5,5] | [-2,2] |
| 8 | [12-20] | 9 | 4 | 1 | 9 | 1 | 1 | [0,20] | [-10,10] | [-2,8] | [-5,15] | [-5,5] | [-2,2] |
| 9 | [12-20] | 9 | 3 | 2 | 8 | 1 | 1 | [0,20] | [-10,10] | [-2,8] | [-5,15] | [-5,5] | [-2,2] |
| 10 | [12-14] | 7 | 1 | 4 | 5 | 2 | 2 | [0,20] | [-10,10] | [-2,8] | [-5,15] | [-5,5] | [-2,2] |
| | | | | | | | | [-3,17] | [-4,6] | [1,7] | [-5,15] | [-1,5] | [0,4] |
| | [15-20] | 9 | 3 | 2 | 8 | 1 | 1 | [0,20] | [-10,10] | [-2,8] | [-5,15] | [-5,5] | [-2,2] |
| | | | | | | | | [-3,17] | [-4,6] | [1,7] | [-5,15] | [-1,5] | [0,4] |
| * 11 | [12-20] | 11 | 4 | 2 | 10 | 1 | 2 | [0,20] | [-10,10] | [-2,8] | [-5,15] | [-5,5] | [-2,2] |
| | | | | | | | | [1,21] | [-1,9] | [-1,5] | [0,20] | [-2,4] | [0,4] |

* These are the only values of r where the best Braunschadel basis betters the best Hofmeister basis for some values of t.

## 3 Discussion

We see that for Hofmeister sets the results for large t (eg t=20) agree with the earlier results given in Table 1 above, with the exception of the r=9 case. The reason for this is that I noticed earlier in 1989 that the cover for the set with $c_{ij}$ = (8,4,1,7,1,0) - given in Table 1 - is "catching up" the cover for the set with $c_{ij}$ = (7,3,3,5,2,2) - given in Table 2: so although the latter is still best for t=20, it will not always remain so.

It turns out that for t<=20, the only cases where the best Braunschadel base has a greater cover than the best Hofmeister base are:

r = 0 for t >= 12
r = 11 for t >= 8

However, examination of the full results suggested that some Braunschadel bases were catching up on their Hofmeister counterparts, and I showed in April 1992, for example, that for r=2, the t=20 Braunschadel basis has a greater cover than the t=20 Hofmeister basis for t>=21.

## 4 Coefficients of powers of t

I returned to look more closely at the "catching up" concept in October 1992, when I came to write up these notes for the first time.

We know that the cover of any Hofmeister or Braunschadel set can be expressed in the regular form as:

$$C = (k_{54}t + c_{54})a_4 + (k_{53}t + c_{53})a_3 + (k_{52}t + c_{52})a_2 + (k_{51}t + c_{51})$$

and, given a particular set for a particular (but reasonably large) value of t, it turns out to be very easy to "guess" what these coefficients are.

As an example, consider the best Braunschadel basis for r=0, t=20:

| | |
|---|---|
| | $k_{ij} = (9, 2, 3, 7, 2, 2)$, $c_{ij} = (2, 2, -1, 3, -1,0)$ |
| gives | $a_i = \{1, 182, 10780, 438441\}$ |
| with | $C = 28491279$ |

We can easily work out that:

$$28491279 = 64a_4 + 39a_3 + 58a_2 + 79$$

from which we can make the educated guess that:

$$C = (3t + 4)a_4 + (2t - 1)a_3 + (3t - 2)a_2 + (4t - 1)$$

thus suggesting coefficients as follows:

| | |
|---|---|
| | $(k_{51}, k_{52}, k_{53}, k_{54}) = (4, 3, 2, 3)$ |
| and | $(c_{51}, c_{52}, c_{53}, c_{54}) = (-1, -2, -1, 4)$ |

We can check out these values by calculating the cover using this formula for other values of t, and checking that these match with the computed values.

The final step is to expand the formula for the cover into a polynomial in t; we find that in this case:

$$C = 162t^4 + 318t^3 + 68t^2 + 4t - 1$$

We can compare this with the corresponding polynomial for the best Hofmeister basis for r=0, which turns out to be:

$$C' = 162t^4 + 312t^3 + 137t^2 + 19t - 2$$

From this it's immediately obvious that $C > C'$ for large enough t, although we know that $C < C'$ for small t.

Some results are given in Tables 4, 5, 6 and 7.

Table 4
Coefficients for the best Hofmeister bases

| r | $c_{21}$ ... $c_{43}$ | | | | | | $k_{51}$-$k_{54}$ | | | | $c_{51}$-$c_{54}$ | | | | coeffs of powers of t | | | | |
|---|---|---|---|---|---|---|---|---|---|---|---|---|---|---|---|---|---|---|---|
| 0 | 2 | 1 | 0 | 1 | 0 | 1 | 5 | 3 | 2 | 3 | -1 | -1 | 0 | 1 | 162 | 303 | 181 | 33 | -1 |
| | 1 | 0 | 0 | 0 | 0 | 0 | 5 | 3 | 2 | 3 | -1 | -1 | -1 | 3 | 162 | 312 | 137 | 19 | -2 |
| 1 | 1 | 0 | 2 | 1 | 1 | 0 | 5 | 3 | 2 | 3 | -1 | 1 | -1 | 2 | 162 | 366 | 252 | 46 | 2 |
| 2 | 2 | 1 | 1 | 1 | 1 | 1 | 5 | 3 | 2 | 3 | -1 | 0 | 0 | 2 | 162 | 411 | 366 | 125 | 11 |
| | 1 | 0 | 2 | 1 | 1 | 0 | 5 | 3 | 2 | 3 | -1 | 1 | -1 | 3 | 162 | 420 | 320 | 68 | 4 |
| 3 | 3 | 1 | 2 | 2 | 1 | 1 | 5 | 3 | 2 | 3 | 0 | 1 | 0 | 2 | 162 | 483 | 504 | 207 | 27 |
| 4 | 3 | 1 | 2 | 2 | 1 | 1 | 5 | 3 | 2 | 3 | 0 | 1 | 0 | 3 | 162 | 537 | 611 | 274 | 39 |
| 5 | 3 | 1 | 2 | 2 | 1 | 1 | 5 | 3 | 2 | 3 | 0 | 1 | 0 | 4 | 162 | 591 | 718 | 341 | 51 |
| 6 | 3 | 1 | 2 | 2 | 1 | 1 | 5 | 3 | 2 | 3 | 0 | 1 | 0 | 5 | 162 | 645 | 825 | 408 | 63 |
| 7 | 7 | 3 | 2 | 5 | 1 | 2 | 5 | 3 | 2 | 3 | 2 | 1 | 1 | 3 | 162 | 690 | 1064 | 700 | 164 |
| | 8 | 4 | 1 | 7 | 1 | 0 | 5 | 3 | 2 | 3 | 2 | 0 | -1 | 7 | 162 | 708 | 886 | 453 | 95 |
| 8 | 7 | 3 | 3 | 5 | 2 | 2 | 5 | 3 | 2 | 3 | 2 | 2 | 1 | 3 | 162 | 744 | 1259 | 918 | 241 |
| | 8 | 4 | 1 | 7 | 1 | 0 | 5 | 3 | 2 | 3 | 2 | 0 | -1 | 8 | 162 | 762 | 978 | 509 | 110 |
| 9 | 7 | 3 | 3 | 5 | 2 | 2 | 5 | 3 | 2 | 3 | 2 | 2 | 1 | 4 | 162 | 798 | 1435 | 1109 | 308 |
| | 8 | 4 | 1 | 7 | 1 | 0 | 5 | 3 | 2 | 3 | 2 | 0 | -1 | 9 | 162 | 816 | 1070 | 565 | 125 |
| 10 | 7 | 3 | 3 | 5 | 2 | 2 | 5 | 3 | 2 | 3 | 2 | 2 | 1 | 5 | 162 | 852 | 1611 | 1300 | 375 |
| | 11 | 6 | 1 | 10 | 1 | 0 | 9 | 3 | 2 | 3 | 8 | 0 | -1 | 9 | 162 | 870 | 1298 | 741 | 180 |
| 11 | 10 | 4 | 3 | 7 | 2 | 2 | 5 | 3 | 2 | 3 | 4 | 2 | 1 | 5 | 162 | 906 | 1851 | 1642 | 533 |
| | 11 | 5 | 2 | 9 | 2 | 0 | 5 | 3 | 2 | 3 | 4 | 1 | -1 | 9 | 162 | 924 | 1565 | 1048 | 267 |

Table 5
Coefficients for the best Braunschadel bases

| r | $c_{21}$ ... $c_{43}$ | | | | | | $k_{51}$-$k_{54}$ | | | | $c_{51}$-$c_{54}$ | | | | coeffs of powers of t | | | | |
|---|---|---|---|---|---|---|---|---|---|---|---|---|---|---|---|---|---|---|---|
| 0 | 2 | 2 | -1 | 3 | -1 | 0 | 4 | 3 | 2 | 3 | -1 | -2 | -1 | 4 | 162 | 318 | 68 | 4 | -1 |
| 1 | 0 | -1 | 2 | -1 | 0 | 1 | 9 | 2 | 2 | 3 | -2 | 1 | 0 | 1 | 162 | 363 | 218 | 35 | -4 |
| | 2 | 2 | -1 | 3 | -1 | 0 | 4 | 3 | 2 | 3 | -1 | -2 | -1 | 5 | 162 | 372 | 84 | 6 | 0 |
| 2 | 0 | -1 | 2 | -1 | 0 | 1 | 9 | 2 | 2 | 3 | -2 | 1 | 0 | 2 | 162 | 417 | 303 | 60 | -6 |
| | 5 | 3 | -1 | 6 | -1 | 0 | 4 | 3 | 2 | 3 | 0 | -2 | -1 | 5 | 162 | 426 | 198 | 12 | -3 |
| 3 | 5 | 3 | 0 | 6 | 0 | 0 | 4 | 3 | 2 | 3 | 0 | -1 | -1 | 5 | 162 | 480 | 363 | 132 | 22 |
| 4 | 3 | 0 | 3 | 2 | 1 | 1 | 9 | 2 | 2 | 3 | 1 | 1 | 0 | 2 | 162 | 525 | 570 | 240 | 32 |
| | 5 | 3 | 0 | 6 | 0 | 0 | 4 | 3 | 2 | 3 | 0 | -1 | -1 | 6 | 162 | 534 | 415 | 155 | 28 |
| 5 | 3 | 0 | 3 | 2 | 1 | 1 | 9 | 2 | 2 | 3 | 1 | 1 | 0 | 3 | 162 | 579 | 691 | 318 | 46 |
| 6 | 3 | 0 | 3 | 2 | 1 | 1 | 9 | 2 | 2 | 3 | 1 | 1 | 0 | 4 | 162 | 633 | 812 | 396 | 60 |
| 7 | 8 | 4 | 1 | 9 | 1 | 0 | 4 | 3 | 2 | 3 | 1 | 0 | -1 | 7 | 162 | 696 | 854 | 460 | 108 |
| 8 | 9 | 4 | 1 | 9 | 1 | 1 | 4 | 3 | 2 | 3 | 2 | 0 | 0 | 6 | 162 | 741 | 1123 | 738 | 188 |
| 9 | 9 | 3 | 2 | 8 | 1 | 1 | 8 | 3 | 2 | 3 | 6 | 1 | 0 | 6 | 162 | 795 | 1324 | 938 | 243 |
| 10 | 7 | 1 | 4 | 5 | 2 | 2 | 7 | 3 | 2 | 3 | 4 | 3 | 1 | 4 | 162 | 840 | 1577 | 1263 | 362 |
| | 9 | 3 | 2 | 8 | 1 | 1 | 8 | 3 | 2 | 3 | 6 | 1 | 0 | 7 | 162 | 849 | 1463 | 1061 | 281 |
| 11 | 11 | 4 | 2 | 10 | 1 | 2 | 4 | 3 | 2 | 3 | 3 | 1 | 1 | 6 | 162 | 912 | 1816 | 1546 | 478 |

## Table 6
### Coefficients for some other good Braunschadel bases for r=11

The sets considered are all those with cover greater than or equal to that of the Braunschadel basis with cij = (7, 1, 4, 5, 2, 2) in the coefficient range ([0,20], [-10,10], [-2,8], [-5,15], [-5,5], [-2,2]).

| $c_{21}$ ... $c_{43}$ | | | | | | $k_{51}$-$k_{54}$ | | | | $c_{51}$-$c_{54}$ | | | | coeffs of powers of t | | | | |
|---|---|---|---|---|---|---|---|---|---|---|---|---|---|---|---|---|---|---|
| 7 | 1 | 4 | 5 | 2 | 2 | 7 | 3 | 2 | 3 | 4 | 3 | 1 | 5 | 162 | 894 | 1767 | 1478 | 439 |
| 9 | 3 | 2 | 8 | 1 | 1 | 8 | 3 | 2 | 3 | 6 | 1 | 0 | 8 | 162 | 903 | 1602 | 1184 | 319 |
| 10 | 3 | 3 | 9 | 1 | 2 | 9 | 2 | 2 | 3 | 8 | 1 | 1 | 5 | 162 | 894 | 1773 | 1518 | 476 |
| 11 | 4 | 2 | 10 | 1 | 2 | 4 | 3 | 2 | 3 | 3 | 1 | 1 | 6 | 162 | 912 | 1816 | 1546 | 478 |
| 11 | 5 | 2 | 12 | 2 | 0 | 4 | 3 | 2 | 3 | 2 | 1 | -1 | 9 | 162 | 912 | 1525 | 1058 | 292 |
| 12 | 5 | 1 | 12 | 1 | 0 | 4 | 3 | 2 | 3 | 3 | 0 | -1 | 10 | 162 | 930 | 1436 | 839 | 226 |
| 12 | 6 | 1 | 13 | 1 | 1 | 4 | 3 | 2 | 3 | 3 | 0 | 0 | 8 | 162 | 903 | 1602 | 1189 | 347 |

## Table 7
### Coefficients for some other good Hofmeister bases for r=0

The sets considered are all those whose cover for t=18 is greater than 18634000 in the (narrow) coefficient range ([-2,2], [-3,1], [0,2], [-6,1],[0,2], [0,1]).

| $c_{21}$ ... $c_{43}$ | | | | | | $k_{51}$-$k_{54}$ | | | | $c_{51}$-$c_{54}$ | | | | coeffs of powers of t | | | | |
|---|---|---|---|---|---|---|---|---|---|---|---|---|---|---|---|---|---|---|
| 1 | 0 | 0 | 0 | 0 | 0 | 5 | 3 | 2 | 3 | -1 | -1 | -1 | 3 | 162 | 312 | 137 | 19 | -2 |
| -2 | -1 | 0 | -2 | 0 | 0 | 5 | 3 | 2 | 3 | -3 | -1 | -1 | 4 | 162 | 312 | 55 | -12 | -8 |
| 2 | 1 | 0 | 1 | 0 | 1 | 5 | 3 | 2 | 3 | -1 | -1 | 0 | 1 | 162 | 303 | 181 | 33 | -1 |
| 0 | 0 | 1 | 0 | 1 | 0 | 5 | 3 | 2 | 3 | -2 | 0 | -1 | 2 | 162 | 294 | 162 | 24 | -2 |
| 0 | -1 | 2 | -1 | 1 | 0 | 6 | 3 | 2 | 3 | -2 | 1 | -1 | 1 | 162 | 294 | 148 | 2 | -2 |
| 0 | 0 | 1 | 0 | 0 | 0 | 6 | 3 | 2 | 3 | -2 | 0 | -1 | 2 | 162 | 294 | 135 | 7 | -2 |
| 0 | 0 | 0 | 0 | 0 | 0 | 5 | 3 | 2 | 3 | -2 | -1 | -1 | 3 | 162 | 294 | 107 | 13 | -2 |
| 1 | 1 | 0 | 0 | 1 | 1 | 5 | 3 | 2 | 3 | -2 | -1 | 0 | 1 | 162 | 285 | 181 | 34 | -1 |
| 1 | 0 | 1 | 0 | 0 | 1 | 8 | 1 | 2 | 3 | 0 | -1 | 1 | 0 | 162 | 285 | 149 | 21 | 0 |
| -2 | -1 | 1 | -2 | 0 | 1 | 6 | 3 | 2 | 3 | -3 | 0 | 0 | 1 | 162 | 285 | 112 | -17 | -8 |
| -2 | -2 | 2 | -3 | 1 | 1 | 5 | 3 | 2 | 3 | -2 | 1 | 0 | 0 | 162 | 285 | 107 | -37 | -4 |
| -2 | -2 | 1 | -3 | 0 | 1 | 5 | 3 | 2 | 3 | -2 | 0 | 0 | 1 | 162 | 285 | 106 | -28 | -9 |
| -2 | -2 | 2 | -2 | 0 | 1 | 5 | 2 | 2 | 3 | -2 | 2 | 0 | 0 | 162 | 285 | 71 | -17 | -6 |
| 2 | 1 | 0 | 1 | 1 | 0 | 5 | 3 | 2 | 3 | -1 | -1 | -1 | 2 | 162 | 276 | 162 | 47 | 2 |
| 2 | 1 | 0 | 1 | 0 | 0 | 8 | 1 | 2 | 3 | 1 | -2 | 0 | 2 | 162 | 276 | 144 | 23 | -1 |
| -1 | -2 | 2 | -3 | 2 | 0 | 5 | 3 | 2 | 3 | -1 | 1 | -1 | 1 | 162 | 276 | 139 | -16 | -3 |
| 2 | 0 | 0 | 0 | 0 | 0 | 8 | 1 | 2 | 3 | 2 | -2 | 0 | 2 | 162 | 276 | 138 | 14 | -2 |
| -1 | -1 | 1 | -2 | 1 | 0 | 6 | 3 | 2 | 3 | -2 | 0 | -1 | 2 | 162 | 276 | 126 | 0 | -6 |
| -1 | -2 | 1 | -3 | 1 | 0 | 5 | 3 | 2 | 3 | -1 | 0 | -1 | 2 | 162 | 276 | 120 | -10 | -6 |
| -1 | -1 | 2 | -1 | 1 | 0 | 6 | 3 | 2 | 3 | -3 | 1 | -1 | 1 | 162 | 276 | 118 | -11 | -3 |
| -1 | -1 | 1 | -1 | 0 | 0 | 5 | 2 | 2 | 3 | -2 | 0 | 0 | 2 | 162 | 276 | 117 | -2 | -4 |
| -1 | -2 | 2 | -2 | 1 | 0 | 5 | 3 | 2 | 3 | -1 | 1 | -1 | 1 | 162 | 276 | 112 | -19 | -1 |
| -1 | 0 | 1 | 0 | 0 | 0 | 6 | 3 | 2 | 3 | -3 | 0 | -1 | 2 | 162 | 276 | 105 | -3 | -2 |
| -1 | 0 | 0 | -1 | 1 | 0 | 5 | 3 | 2 | 3 | -3 | -1 | -1 | 3 | 162 | 276 | 104 | 28 | -8 |
| -1 | -1 | 2 | 0 | 0 | 0 | 5 | 2 | 2 | 3 | -2 | 2 | -1 | 1 | 162 | 276 | 82 | -5 | -1 |
| -1 | -1 | 2 | -1 | 0 | 0 | 5 | 2 | 2 | 3 | -2 | 2 | -1 | 1 | 162 | 276 | 82 | -8 | -2 |
| -1 | -1 | 0 | -2 | 0 | 0 | 8 | 1 | 2 | 3 | -1 | -2 | 0 | 3 | 162 | 276 | 80 | -10 | -5 |
| -1 | 0 | 0 | -1 | 0 | 0 | 6 | 3 | 2 | 3 | -2 | -1 | -1 | 3 | 162 | 276 | 77 | 5 | -4 |
| -1 | -2 | 0 | -3 | 0 | 0 | 5 | 3 | 2 | 3 | -1 | -1 | -1 | 3 | 162 | 276 | 65 | -18 | -7 |

From these tables, we see that the best Braunschadel basis for t=20 will eventually better the best Hofmeister basis for t=20 only in the case of r = 0, 1 or 2. Of course, this does *not* mean that Hofmeister is always best in all other cases - simply that we have not yet found a Braunschadel basis with better coefficients.

The case of r=11 is particularly interesting, since we know that a Braunschadel basis is best for t = 8 to 20, whereas Tables 4 and 5 show that the best Hofmeister basis for t=20 must eventually dominate. But Table 6 reveals another Braunschadel basis with yet a higher coefficient of $t^3$, and we find that the optimal basis for r=11 is likely to switch from Hofmeister to Braunschadel and back again as follows:

| $t_{range}$ | type | $c_{21}$ - $c_{43}$ | coeffs of powers of t | | | | |
|---|---|---|---|---|---|---|---|
| [2-7] | H | (10, 4, 3, 7, 2, 2) | 162 | 906 | 1851 | 1642 | 533 |
| [8-22] | B | (11, 4, 2, 10, 1, 2) | 162 | 912 | 1816 | 1546 | 748 |
| 23 | H | (11, 5, 2, 9, 2, 0) | 162 | 924 | 1565 | 1048 | 267 |
| [24-...] | B | (12, 5, 1, 12, 1, 0) | 162 | 930 | 1436 | 839 | 226 |

[ Later results confirm this behaviour - see section 6 and Appendix 1 ]

According to Mossige's thesis (see Mossige, Svein, [6], p38), the formula for the cover of any Hofmeister basis must conform to one of seven possibilities $n_i(x,P)$. The coefficients $c_{5j}$ in these possibilities depend on the coefficients $c_{ij}$ for i<=4, but the coefficients $k_{5j}$ are explicit, as follows:

| | |
|---|---|
| $n_1$ | (8, 1, 0, 3) |
| $n_2$ | (5, 3, 2, 3) |
| $n_3$ | (9, 3, 2, 3) |
| $n_4$ | (9, 3, 2, 3) |
| $n_5$ | (8, 1, 2, 3) |
| $n_6$ | (9, 1, 2, 3) |
| $n_7$ | (5, 2, 2, 3) |

Examining Tables 4 and 7 above show that the bases in these tables conform to these requirements with the exception of those with coefficients (6, 3, 2, 3). This was at first a mystery - even after I had taken into account Selmer's correction to Mossige's thesis (Selmer, E.S., [8]; see also Appendix 3) - but closer examination of my own notes after first reading the thesis in March 1992 (see particularly Mossige, Svein, [6] pp43-45) show that a further case omitted by Mossige and not noticed by Selmer can generate coefficients (6, 3, 2, 3).

## 5 The future

Future work should really be based on "Proposition 6.1" (and its analogue for Braunschadel bases - probably very similar) from Mossige's thesis: this gives conditions which, if met by a Hofmeister basis, determine that it is admissible, and also give a formula for its cover. It is possible (although by no means certain) that using these conditions would speed up exp21 so that wider ranges of $c_{ij}$ might be investigated (or, at least, easily rejected), and larger values of t investigated.

Mossige's 2.008... result shows that there is no "finite set of formulae" for d=4, and I now conjecture instead that all maximal sets for sufficiently large h are either Hofmeister or Braunschadel - and that the coefficient sets $c_{ij}$ will, for a fixed value of r, change from time to time, slowly but surely moving away from (0,0,0,0,0,0) without limit. Some interesting further questions can also be posed:

Will one or the other form eventually dominate?
Will one or the other form eventually dominate for a given value of r?
What is the "sufficiently large value" of h?

## 6 Afterword - further results

Following a request from Selmer in December 1992 for "good" bases for larger values of t, I looked again at the program I had been using to investigate these bases. I was able to incorporate some worthwhile improvements* to the algorithms used, allowing me to extend the searches for good bases to around t=58. Beyond this point the cover exceeds $2^{31}$, and a major coding change to double-length working will be needed to progress further.

[ *The first improvement was to check admissibility criteria: if a candidate set $\{1, a_2, a_3, a_4\}$ proves to be inadmissible, then all other candidate sets $\{1, a_2, a_3, a_4'\}$ are rejected without further ado.
The second improvement was to the code which checks the cover of a "good" set. Advantage is taken of the fact that if x has a generation:

$$x = c_4a_4 + c_3a_3 + c_2a_2 + c_1$$

then so do all $(x-ia_4)$ for $0 < i <= c_4$.
This significantly reduces the time taken to determine the full cover of a good set. ]

These new results are summarised in the following tables*, which should be considered as extensions of tables 2, 3, 4 and 5 above.

I have also calculated and compared the covers of the best Hofmeister and Braunschadel bases for each value of t and r, and the results are summarised in Appendix 1 which is, effectively, a table of the best bases found so far for k=4 and 144 <= h <= 707.

[ *These results have now (March 93) been further confirmed by extended searches.
For Hofmeister bases, the searches covered the range:
$c_{21}$: [ -8, 21]  $c_{31}$: [-11, 18]  $c_{32}$: [-15, 14]  $c_{41}$: [ -8, 21]  $c_{42}$: [-14, 15]  $c_{43}$: [-16, 13]
and for Braunschadel bases:
$c_{21}$: [ -7, 22]  $c_{31}$: [-10, 19]  $c_{32}$: [-15, 14]  $c_{41}$: [ -5, 24]  $c_{42}$: [-14, 15]  $c_{43}$: [-16, 13]
Any basis whose cover for t=58 exceeded that of the best t=21 basis was printed out.
No improved bases for any t in the range 21 <= t <= 58 were found. ]

Table 8

Probably the best Hofmeister bases for 21<=t<=58

| | | best coefficients | | | | | | range of coefficients checked | | | | | |
|---|---|---|---|---|---|---|---|---|---|---|---|---|---|
| r | $t_{range}$ | $c_{21}$ | $c_{31}$ | $c_{32}$ | $c_{41}$ | $c_{42}$ | $c_{43}$ | $c_{21}$ | $c_{31}$ | $c_{32}$ | $c_{41}$ | $c_{42}$ | $c_{43}$ |
| 0 | [21-22] | 1 | 0 | 0 | 0 | 0 | 0 | [-2,9] | [-1,6] | [-6,1] | [-2,10] | [-2,1] | [-6,1] |
| | [23-39] | 2 | 1 | -1 | 2 | 0 | -2 | | | | | | |
| | [40-58] | 6 | 5 | -4 | 8 | -1 | -5 | | | | | | |
| 1 | [21-30] | 1 | 0 | 2 | 1 | 1 | 0 | [-2,9] | [-1,6] | [-4,3] | [-2,10] | [-1,2] | [-6,1] |
| | [31-34] | 2 | 1 | -1 | 2 | 0 | -2 | | | | | | |
| | [35-58] | 6 | 4 | -3 | 7 | 0 | -5 | | | | | | |
| 2 | [21-28] | 1 | 0 | 2 | 1 | 1 | 0 | [-2,9] | [-1,6] | [-4,3] | [-2,10] | [-1,2] | [-6,1] |
| | [29-40] | 5 | 3 | -1 | 5 | 0 | -2 | | | | | | |
| | [41-58] | 6 | 4 | -3 | 7 | 0 | -5 | | | | | | |
| 3 | [21-48] | 3 | 1 | 2 | 2 | 1 | 1 | [-1,10] | [-1,6] | [-4,3] | [-1,11] | [-1,2] | [-6,1] |
| | [49-58] | 6 | 4 | -3 | 7 | 0 | -5 | | | | | | |
| 4 | [21-48] | 3 | 1 | 2 | 2 | 1 | 1 | [1,12] | [0,7] | [-4,3] | [0,12] | [-1,2] | [-6,1] |
| | [49-58] | 9 | 6 | -3 | 10 | 0 | -5 | | | | | | |
| 5 | [21-45] | 3 | 1 | 2 | 2 | 1 | 1 | [1,12] | [0,7] | [-4,3] | [0,12] | [0,3] | [-6,1] |
| | [46-58] | 9 | 5 | -2 | 9 | 1 | -5 | | | | | | |
| 6 | [21-49] | 3 | 1 | 2 | 2 | 1 | 1 | [1,12] | [0,7] | [-4,3] | [0,12] | [0,3] | [-6,1] |
| | [50-58] | 9 | 5 | -2 | 9 | 1 | -5 | | | | | | |
| 7 | [21-36] | 8 | 4 | 1 | 7 | 1 | 0 | [3,14] | [1,8] | [-3,4] | [2,14] | [0,3] | [-5,2] |
| | [37-58] | 9 | 5 | 0 | 9 | 1 | -2 | | | | | | |
| 8 | [21-40] | 8 | 4 | 1 | 7 | 1 | 0 | [3,14] | [1,8] | [-3,4] | [2,14] | [0,3] | [-5,2] |
| | [41-58] | 9 | 5 | 0 | 9 | 1 | -2 | | | | | | |
| 9 | [21] | 7 | 3 | 3 | 5 | 2 | 2 | [4,15] | [2,9] | [-2,5] | [3,15] | [0,3] | [-4,3] |
| | [22-29] | 8 | 4 | 1 | 7 | 1 | 0 | | | | | | |
| | [30-58] | 12 | 7 | -1 | 12 | 1 | -3 | | | | | | |
| 10 | [21-34] | 11 | 6 | 1 | 10 | 1 | 0 | [6,17] | [3,10] | [-3,4] | [5,17] | [0,3] | [-5,2] |
| | [35-58] | 12 | 7 | -1 | 12 | 1 | -3 | | | | | | |
| 11 | [21-40] | 11 | 5 | 2 | 9 | 2 | 0 | [6,17] | [3,10] | [-3,4] | [5,17] | [0,3] | [-5,2] |
| | [41-58] | 12 | 7 | -1 | 12 | 1 | -3 | | | | | | |

Notes: The ranges given apply to all values of t for the given value of r. Although not as extensive as the ranges used when checking values of t<=20, some limited experiments with much wider ranges for certain values of r suggest that they are adequate [but see note on previous page: much more extensive searches now give greater confidence in these results (March 93)].

Table 9
Probably the best Braunschadel bases for 21<=t<=58

| | | best coefficients | | | | | | range(s) of coefficients checked | | | | | |
|---|---|---|---|---|---|---|---|---|---|---|---|---|---|
| r | $t_{range}$ | $c_{21}$ | $c_{31}$ | $c_{32}$ | $c_{41}$ | $c_{42}$ | $c_{43}$ | $c_{21}$ | $c_{31}$ | $c_{32}$ | $c_{41}$ | $c_{42}$ | $c_{43}$ |
| 0 | [21-27] | 2 | 2 | -1 | 3 | -1 | 0 | [-3,10] | [-1,7] | [-6,2] | [-2,12] | [-3,1] | [-6,2] |
| | [28-53] | 3 | 3 | -2 | 5 | -1 | -2 | | | | | | |
| | [54-58] | 4 | 4 | -3 | 7 | -1 | -4 | | | | | | |
| 1 | [21-31] | 2 | 2 | -1 | 3 | -1 | 0 | [-2,11] | [0,8] | [-7,1] | [0,14] | [-3,1] | [-7,1] |
| | [32-37] | 3 | 3 | -2 | 5 | -1 | -2 | | | | | | |
| | [38-58] | 7 | 6 | -4 | 11 | -1 | -5 | | | | | | |
| 2 | [21-31] | 5 | 3 | -1 | 6 | -1 | 0 | [0,13] | [1,9] | [-7,1] | [2,16] | [-3,1] | [-7,1] |
| | [32-43] | 6 | 4 | -2 | 8 | -1 | -2 | | | | | | |
| | [44-58] | 7 | 6 | -4 | 11 | -1 | -5 | | | | | | |
| 3 | [21-30] | 5 | 3 | 0 | 6 | 0 | 0 | [0,13] | [1,9] | [-6,2] | [2,16] | [-3,1] | [-7,1] |
| | [31-50] | 6 | 4 | -1 | 8 | 0 | -2 | | | | | | |
| | [51-58] | 7 | 6 | -4 | 11 | -1 | -5 | | | | | | |
| 4 | [21-34] | 5 | 3 | 0 | 6 | 0 | 0 | [1,14] | [1,9] | [-6,2] | [3,17] | [-3,1] | [-7,1] |
| | [35-47] | 6 | 4 | -1 | 8 | 0 | -2 | | | | | | |
| | [48-58] | 10 | 7 | -4 | 14 | -1 | -5 | | | | | | |
| 5 | [21-25] | 3 | 0 | 3 | 2 | 1 | 1 | [0,13] | [0,8] | [-4,4] | [1,15] | [-2,2] | [-6,2] |
| | [26-37] | 5 | 3 | 0 | 6 | 0 | 0 | | | | | | |
| | [38-41] | 6 | 4 | -1 | 8 | 0 | -2 | | | | | | |
| | [42-58] | 10 | 7 | -3 | 14 | 0 | -5 | | | | | | |
| 6 | [21-38] | 8 | 4 | 0 | 9 | 0 | 0 | [3,16] | [2,10] | [-6,2] | [5,19] | [-2,2] | [-7,1] |
| | [39-47] | 9 | 5 | -1 | 11 | 0 | -2 | | | | | | |
| | [48-58] | 10 | 7 | -3 | 14 | 0 | -5 | | | | | | |
| 7 | [21-37] | 8 | 4 | 1 | 9 | 1 | 0 | [3,16] | [2,10] | [-5,3] | [5,19] | [-2,2] | [-7,1] |
| | [38-55] | 9 | 5 | 0 | 11 | 1 | -2 | | | | | | |
| | [56-58] | 10 | 7 | -3 | 14 | 0 | -5 | | | | | | |
| 8 | [21] | 9 | 4 | 1 | 9 | 1 | 1 | [5,18] | [2,10] | [-5,3] | [6,20] | [-2,2] | [-6,2] |
| | [22-40] | 8 | 4 | 1 | 9 | 1 | 0 | | | | | | |
| | [41-52] | 9 | 5 | 0 | 11 | 1 | -2 | | | | | | |
| | [53-58] | 13 | 8 | -3 | 17 | 0 | -5 | | | | | | |
| 9 | [21-27] | 9 | 3 | 2 | 8 | 1 | 1 | [5,18] | [2,10] | [-4,4] | [6,20] | [-1,3] | [-6,2] |
| | [28-49] | 10 | 4 | 1 | 10 | 1 | -1 | | | | | | |
| | [50-58] | 13 | 8 | -2 | 17 | 1 | -5 | | | | | | |
| 10 | [21-24] | 9 | 3 | 2 | 8 | 1 | 1 | [5,18] | [2,10] | [-4,4] | [6,20] | [-1,3] | [-6,2] |
| | [25-37] | 11 | 5 | 1 | 12 | 1 | 0 | | | | | | |
| | [38-51] | 10 | 4 | 1 | 10 | 1 | -1 | | | | | | |
| | [52-58] | 13 | 8 | -2 | 17 | 1 | -5 | | | | | | |
| 11 | [21-22] | 11 | 4 | 2 | 10 | 1 | 2 | [6,19] | [1,9] | [-3,5] | [5,19] | [-1,3] | [-4,4] |
| | [23-48] | 12 | 5 | 1 | 12 | 1 | 0 | | | | | | |
| | [49-58] | 13 | 6 | 0 | 14 | 1 | -2 | | | | | | |

Note: The ranges given apply to all values of t for the given value of r.

Table 10

Coefficients for the best Hofmeister bases for 21<=t<=58

| r | $c_{21}$ ... $c_{43}$ | | | | | | $k_{51}$-$k_{54}$ | | | | $c_{51}$-$c_{54}$ | | | | coeffs of powers of t | | | | |
|---|---|---|---|---|---|---|---|---|---|---|---|---|---|---|---|---|---|---|---|
| 0 | 1 | 0 | 0 | 0 | 0 | 0 | 5 | 3 | 2 | 3 | -1 | -1 | -1 | 3 | 162 | 312 | 137 | 19 | -2 |
| | 2 | 1 | -1 | 2 | 0 | -2 | 5 | 3 | 2 | 3 | -1 | -2 | -3 | 7 | 162 | 330 | -269 | 49 | 26 |
| | 6 | 5 | -4 | 8 | -1 | -5 | 9 | 3 | 2 | 3 | 1 | -5 | -6 | 14 | 162 | 375 | -2067 | 907 | 1443 |
| 1 | 1 | 0 | 2 | 1 | 1 | 0 | 5 | 3 | 2 | 3 | -1 | 1 | -1 | 2 | 162 | 366 | 252 | 46 | 2 |
| | 2 | 1 | -1 | 2 | 0 | -2 | 5 | 3 | 2 | 3 | -1 | -2 | -3 | 8 | 162 | 384 | -303 | 56 | 30 |
| | 6 | 4 | -3 | 7 | 0 | -5 | 5 | 3 | 2 | 3 | 0 | -4 | -6 | 14 | 162 | 429 | -1875 | 444 | 1138 |
| 2 | 1 | 0 | 2 | 1 | 1 | 0 | 5 | 3 | 2 | 3 | -1 | 1 | -1 | 3 | 162 | 420 | 320 | 68 | 4 |
| | 5 | 3 | -1 | 5 | 0 | -2 | 9 | 3 | 2 | 3 | 2 | -2 | -3 | 8 | 162 | 438 | -183 | -57 | 70 |
| | 6 | 4 | -3 | 7 | 0 | -5 | 5 | 3 | 2 | 3 | 0 | -4 | -6 | 15 | 162 | 483 | -2002 | 460 | 1215 |
| 3 | 3 | 1 | 2 | 2 | 1 | 1 | 5 | 3 | 2 | 3 | 0 | 1 | 0 | 2 | 162 | 483 | 504 | 207 | 27 |
| | 6 | 4 | -3 | 7 | 0 | -5 | 5 | 3 | 2 | 3 | 0 | -4 | -6 | 16 | 162 | 537 | -2129 | 476 | 1292 |
| 4 | 3 | 1 | 2 | 2 | 1 | 1 | 5 | 3 | 2 | 3 | 0 | 1 | 0 | 3 | 162 | 537 | 611 | 274 | 39 |
| | 9 | 6 | -3 | 10 | 0 | -5 | 9 | 3 | 2 | 3 | 4 | -4 | -6 | 16 | 162 | 591 | -1982 | -313 | 1934 |
| 5 | 3 | 1 | 2 | 2 | 1 | 1 | 5 | 3 | 2 | 3 | 0 | 1 | 0 | 4 | 162 | 591 | 718 | 341 | 51 |
| | 9 | 5 | -2 | 9 | 1 | -5 | 5 | 3 | 2 | 3 | 2 | -3 | -6 | 16 | 162 | 645 | -1736 | -762 | 1381 |
| 6 | 3 | 1 | 2 | 2 | 1 | 1 | 5 | 3 | 2 | 3 | 0 | 1 | 0 | 5 | 162 | 645 | 825 | 408 | 63 |
| | 9 | 5 | -2 | 9 | 1 | -5 | 5 | 3 | 2 | 3 | 2 | -3 | -6 | 17 | 162 | 699 | -1827 | -819 | 1464 |
| 7 | 8 | 4 | 1 | 7 | 1 | 0 | 5 | 3 | 2 | 3 | 2 | 0 | -1 | 7 | 162 | 708 | 886 | 453 | 95 |
| | 9 | 5 | 0 | 9 | 1 | -2 | 5 | 3 | 2 | 3 | 2 | -1 | -3 | 11 | 162 | 726 | 240 | -234 | 66 |
| 8 | 8 | 4 | 1 | 7 | 1 | 0 | 5 | 3 | 2 | 3 | 2 | 0 | -1 | 8 | 162 | 762 | 978 | 509 | 110 |
| | 9 | 5 | 0 | 9 | 1 | -2 | 5 | 3 | 2 | 3 | 2 | -1 | -3 | 12 | 162 | 780 | 266 | -252 | 74 |
| 9 | 7 | 3 | 3 | 5 | 2 | 2 | 5 | 3 | 2 | 3 | 2 | 2 | 1 | 4 | 162 | 798 | 1435 | 1109 | 308 |
| | 8 | 4 | 1 | 7 | 1 | 0 | 5 | 3 | 2 | 3 | 2 | 0 | -1 | 9 | 162 | 816 | 1070 | 565 | 125 |
| | 12 | 7 | -1 | 12 | 1 | -3 | 5 | 3 | 2 | 3 | 3 | -2 | -4 | 15 | 162 | 861 | -223 | -939 | 584 |
| 10 | 11 | 6 | 1 | 10 | 1 | 0 | 9 | 3 | 2 | 3 | 8 | 0 | -1 | 9 | 162 | 870 | 1298 | 741 | 180 |
| | 12 | 7 | -1 | 12 | 1 | -3 | 5 | 3 | 2 | 3 | 3 | -2 | -4 | 16 | 162 | 915 | -224 | -1002 | 623 |
| 11 | 11 | 5 | 2 | 9 | 2 | 0 | 5 | 3 | 2 | 3 | 4 | 1 | -1 | 9 | 162 | 924 | 1565 | 1048 | 267 |
| | 12 | 7 | -1 | 12 | 1 | -3 | 5 | 3 | 2 | 3 | 3 | -2 | -4 | 17 | 162 | 969 | -225 | -1065 | 662 |

# Table 11

Coefficients for the best Braunschadel bases for 21<=t<=58

| r | | $c_{21}$ ... $c_{43}$ | | | | | $k_{51}$-$k_{54}$ | | | | $c_{51}$-$c_{54}$ | | | | coeffs of powers of t | | | | |
|---|---|---|---|---|---|---|---|---|---|---|---|---|---|---|---|---|---|---|---|
| 0 | 2 | 2 | -1 | 3 | -1 | 0 | 4 | 3 | 2 | 3 | -1 | -2 | -1 | 4 | 162 | 318 | 68 | 4 | -1 |
| | 3 | 3 | -2 | 5 | -1 | -2 | 4 | 3 | 2 | 3 | -1 | -3 | -3 | 8 | 162 | 336 | -432 | 121 | 63 |
| | 4 | 4 | -3 | 7 | -1 | -4 | 8 | 3 | 2 | 3 | 0 | -4 | -5 | 12 | 162 | 354 | -1400 | 642 | 444 |
| 1 | 2 | 2 | -1 | 3 | -1 | 0 | 4 | 3 | 2 | 3 | -1 | -2 | -1 | 5 | 162 | 372 | 84 | 6 | 0 |
| | 3 | 3 | -2 | 5 | -1 | -2 | 4 | 3 | 2 | 3 | -1 | -3 | -3 | 9 | 162 | 390 | -482 | 133 | 71 |
| | 7 | 6 | -4 | 11 | -1 | -5 | 4 | 3 | 2 | 3 | 0 | -5 | -6 | 15 | 162 | 435 | -2207 | 851 | 1807 |
| 2 | 5 | 3 | -1 | 6 | -1 | 0 | 4 | 3 | 2 | 3 | 0 | -2 | -1 | 5 | 162 | 426 | 198 | 12 | -3 |
| | 6 | 4 | -2 | 8 | -1 | -2 | 4 | 3 | 2 | 3 | 0 | -3 | -3 | 9 | 162 | 444 | -368 | -63 | 168 |
| | 7 | 6 | -4 | 11 | -1 | -5 | 4 | 3 | 2 | 3 | 0 | -5 | -6 | 16 | 162 | 489 | -2350 | 884 | 1921 |
| 3 | 5 | 3 | 0 | 6 | 0 | 0 | 4 | 3 | 2 | 3 | 0 | -1 | -1 | 5 | 162 | 480 | 363 | 132 | 22 |
| | 6 | 4 | -1 | 8 | 0 | -2 | 4 | 3 | 2 | 3 | 0 | -2 | -3 | 9 | 162 | 498 | -179 | -60 | 102 |
| | 7 | 6 | -4 | 11 | -1 | -5 | 4 | 3 | 2 | 3 | 0 | -5 | -6 | 17 | 162 | 543 | -2493 | 917 | 2035 |
| 4 | 5 | 3 | 0 | 6 | 0 | 0 | 4 | 3 | 2 | 3 | 0 | -1 | -1 | 6 | 162 | 534 | 415 | 155 | 28 |
| | 6 | 4 | -1 | 8 | 0 | -2 | 4 | 3 | 2 | 3 | 0 | -2 | -3 | 10 | 162 | 552 | -193 | -67 | 114 |
| | 10 | 7 | -4 | 14 | -1 | -5 | 4 | 3 | 2 | 3 | 1 | -5 | -6 | 17 | 162 | 597 | -2352 | -22 | 3022 |
| 5 | 3 | 0 | 3 | 2 | 1 | 1 | 9 | 2 | 2 | 3 | 1 | 1 | 0 | 3 | 162 | 579 | 691 | 318 | 46 |
| | 5 | 3 | 0 | 6 | 0 | 0 | 4 | 3 | 2 | 3 | 0 | -1 | -1 | 7 | 162 | 588 | 467 | 178 | 34 |
| | 6 | 4 | -1 | 8 | 0 | -2 | 4 | 3 | 2 | 3 | 0 | -2 | -3 | 11 | 162 | 606 | -207 | -74 | 126 |
| | 10 | 7 | -3 | 14 | 0 | -5 | 4 | 3 | 2 | 3 | 1 | -4 | -6 | 17 | 162 | 651 | -2076 | -439 | 2292 |
| 6 | 8 | 4 | 0 | 9 | 0 | 0 | 4 | 3 | 2 | 3 | 1 | -1 | -1 | 7 | 162 | 642 | 635 | 245 | 52 |
| | 9 | 5 | -1 | 11 | 0 | -2 | 4 | 3 | 2 | 3 | 1 | -2 | -3 | 11 | 162 | 660 | -39 | -251 | 204 |
| | 10 | 7 | -3 | 14 | 0 | -5 | 4 | 3 | 2 | 3 | 1 | -4 | -6 | 18 | 162 | 705 | -2183 | -483 | 2421 |
| 7 | 8 | 4 | 1 | 9 | 1 | 0 | 4 | 3 | 2 | 3 | 1 | 0 | -1 | 7 | 162 | 696 | 854 | 460 | 108 |
| | 9 | 5 | 0 | 11 | 1 | -2 | 4 | 3 | 2 | 3 | 1 | -1 | -3 | 11 | 162 | 714 | 204 | -179 | 87 |
| | 10 | 7 | -3 | 14 | 0 | -5 | 4 | 3 | 2 | 3 | 1 | -4 | -6 | 19 | 162 | 759 | -2290 | -527 | 2550 |
| 8 | 9 | 4 | 1 | 9 | 1 | 1 | 4 | 3 | 2 | 3 | 2 | 0 | 0 | 6 | 162 | 741 | 1123 | 738 | 188 |
| | 8 | 4 | 1 | 9 | 1 | 0 | 4 | 3 | 2 | 3 | 1 | 0 | -1 | 8 | 162 | 750 | 942 | 516 | 125 |
| | 9 | 5 | 0 | 11 | 1 | -2 | 4 | 3 | 2 | 3 | 1 | -1 | -3 | 12 | 162 | 768 | 226 | -193 | 97 |
| | 13 | 8 | -3 | 17 | 0 | -5 | 4 | 3 | 2 | 3 | 2 | -4 | -6 | 19 | 162 | 813 | -2095 | -1504 | 3404 |
| 9 | 9 | 3 | 2 | 8 | 1 | 1 | 8 | 3 | 2 | 3 | 6 | 1 | 0 | 6 | 162 | 795 | 1324 | 938 | 243 |
| | 10 | 4 | 1 | 10 | 1 | -1 | 8 | 3 | 2 | 3 | 6 | 0 | -2 | 10 | 162 | 813 | 854 | 232 | 38 |
| | 13 | 8 | -2 | 17 | 1 | -5 | 4 | 3 | 2 | 3 | 2 | -3 | -6 | 19 | 162 | 867 | -1765 | -1869 | 2351 |
| 10 | 9 | 3 | 2 | 8 | 1 | 1 | 8 | 3 | 2 | 3 | 6 | 1 | 0 | 7 | 162 | 849 | 1463 | 1061 | 281 |
| | 11 | 5 | 1 | 12 | 1 | 0 | 4 | 3 | 2 | 3 | 2 | 0 | -1 | 9 | 162 | 858 | 1252 | 724 | 193 |
| | 10 | 4 | 1 | 10 | 1 | -1 | 8 | 3 | 2 | 3 | 6 | 0 | -2 | 11 | 162 | 867 | 927 | 255 | 44 |
| | 13 | 8 | -2 | 17 | 1 | -5 | 4 | 3 | 2 | 3 | 2 | -3 | -6 | 20 | 162 | 921 | -1836 | -1978 | 2471 |
| 11 | 11 | 4 | 2 | 10 | 1 | 2 | 4 | 3 | 2 | 3 | 3 | 1 | 1 | 6 | 162 | 912 | 1816 | 1546 | 478 |
| | 12 | 5 | 1 | 12 | 1 | 0 | 4 | 3 | 2 | 3 | 3 | 0 | -1 | 10 | 162 | 930 | 1436 | 839 | 226 |
| | 13 | 6 | 0 | 14 | 1 | -2 | 8 | 3 | 2 | 3 | 8 | -1 | -3 | 14 | 162 | 948 | 588 | -420 | 187 |

## 7 Sequences of bases

Examination of these results shows that from time to time *sequences* of good bases appear, whose coefficients differ one from the other by a constant amount. One sequence, particularly prominent in the Braunschadel results, has a coefficient difference set of (1, 1, -1, 2, 0, -2); examples include:

| | | | | | | | |
|---|---|---|---|---|---|---|---|
| 0 | [21-27] | 2 | 2 | -1 | 3 | -1 | 0 |
| | [28-53] | 3 | 3 | -2 | 5 | -1 | -2 |
| | [54-58] | 4 | 4 | -3 | 7 | -1 | -4 |
| 7 | [21-37] | 8 | 4 | 1 | 9 | 1 | 0 |
| | [38-55] | 9 | 5 | 0 | 11 | 1 | -2 |
| 11 | [21-22] | 11 | 4 | 2 | 10 | 1 | 2 |
| | [23-48] | 12 | 5 | 1 | 12 | 1 | 0 |
| | [49-58] | 13 | 6 | 0 | 14 | 1 | -2 |

and the sequence also appears in the Hofmeister results:

| | | | | | | | |
|---|---|---|---|---|---|---|---|
| 0 | [21-22] | 1 | 0 | 0 | 0 | 0 | 0 |
| | [23-39] | 2 | 1 | -1 | 2 | 0 | -2 |
| 7 | [21-36] | 8 | 4 | 1 | 7 | 1 | 0 |
| | [37-58] | 9 | 5 | 0 | 9 | 1 | -2 |

In his thesis, Mossige gives criteria which the coefficients of a Hofmeister basis must satisfy in order to form a basis, and also a formula for the h-range when these criteria are satisfied. It is interesting to apply these criteria and formulae to the sequence of Hofmeister bases defined by the parametrised coefficient set:

$$(8+p, 4+p, 1-p, 7+2p, 1, -2p) \qquad \text{for} \qquad h = 12t + 7$$

Using the notation of (Mossige, Svein, [6] page 37) we have:

$P = (d_1, d_2, d_3, d_4, d_5, d_6)$ where

$d_1 = 8+p$
$d_2 = 4+p$
$d_3 = 1-p$
$d_4 = 7+2p$
$d_5 = 1$
$d_6 = -2p$

We first check that both the "admissibility" conditions are always >= 0:

$f_8(P) = 0$
$f_9(P) = 3p+4$

Next, we check that there is no need to take Selmer's correction into account:

$-d_1-3d_2+3d_4 = 2p+1$
$-2d_3+3d_5+3 = 2p+4$

and now evaluate $f_i(x, P)$ for i=1 to 7:

$f_1(x,P) >= 0 => x <= 4p+7$ [this is $f_{1A}$, since the condition for $f_{1B}$ evaluates to $2p+3<0$]
$f_2(x,P) >= 0 => x <= 4p+6$
$f_3(x,P) >= 0 => x <= 5p+6$
$f_4(x,P) >= 0 => x <= 5p+8$
$f_5(x,P) >= 0 => x <= 5p+7$
$f_6(x,P) >= 0 => x <= 4p+6$
$f_7(x,P) >= 0 => x <= 4p+6$

Clearly, the maximum permissible value of x is 4p+6, and so $n_2$, $n_6$ and $n_7$ are candidate formulae for the cover. Of these, $n_2$ has the lowest coefficient of $a_3$ and so defines the value of the h-range:

$$n_2(4p+6,P) = (3t + (4p+7))a_4 + (2t - (2p+1))a_3 + (3t - p)a_2 + (5t+ 2)$$

Somewhat reassuringly, this means that:

$(k_{51},k_{52},k_{53},k_{54}) = (5,3,2,3)$
and $(c_{51},c_{52},c_{53},c_{54}) = (2, -p, -(2p+1), (4p+7))$

which is confirmed by the appropriate entries in Table 6 for r=7.

Finally, we can calculate the coefficient X of $t^3$ in the h-range; using the formula of Appendix 2 we have:

$$X = (k_{54}k_{43}k_{32}c_{21} + k_{54}k_{43}c_{32}k_{21} + k_{54}c_{43}k_{32}k_{21} + c_{54}k_{43}k_{32}k_{21} + k_{54}k_{43}k_{31} + k_{54}k_{42}k_{21} + k_{53}k_{32}k_{21})$$
$$= 18(8+p) + 54(1-p) - 162p + 54(4p+7)$$
$$= 18p + 708$$

This shows that given any value X we can find a Hofmeister basis (valid for large enough t) whose cover is given by:

$$162t^4 + X't^3 + O(t^2) \quad \text{for some } X' >= X$$

In other words, there is no limit to the value of the coefficient of $t^3$ in the formula for the cover.

**Appendix 1**

This table summarises the best bases found so far for 144 <= h <= 707.

Probably the best Braunschadel or Hofmeister bases for 12<=t<=58

| r | $t_{range}$ | type | $c_{21}$ | $c_{31}$ | $c_{32}$ | $c_{41}$ | $c_{42}$ | $c_{43}$ |
|---|---|---|---|---|---|---|---|---|
| 0 | [12-27] | B | 2 | 2 | -1 | 3 | -1 | 0 |
| | [28-41] | B | 3 | 3 | -2 | 5 | -1 | -2 |
| | [42-58] | H | 6 | 5 | -4 | 8 | -1 | -5 |
| 1 | [12-28] | H | 1 | 0 | 2 | 1 | 1 | 0 |
| | [29-31] | B | 2 | 2 | -1 | 3 | -1 | 0 |
| | [32-35] | B | 3 | 3 | -2 | 5 | -1 | -2 |
| | [36-54] | H | 6 | 4 | -3 | 7 | 0 | -5 |
| | [55-58] | B | 7 | 6 | -4 | 11 | -1 | -5 |
| 2 | [12-20] | H | 1 | 0 | 2 | 1 | 1 | 0 |
| | [21-31] | B | 5 | 3 | -1 | 6 | -1 | 0 |
| | [32-41] | B | 6 | 4 | -2 | 8 | -1 | -2 |
| | [42-56] | H | 6 | 4 | -3 | 7 | 0 | -5 |
| | [57-58] | B | 7 | 6 | -4 | 11 | -1 | -5 |
| 3 | [12-45] | H | 3 | 1 | 2 | 2 | 1 | 1 |
| | [46-49] | B | 6 | 4 | -1 | 8 | 0 | -2 |
| | [50-58] | H | 6 | 4 | -3 | 7 | 0 | -5 |
| 4 | [12-48] | H | 3 | 1 | 2 | 2 | 1 | 1 |
| | [49-58] | H | 9 | 6 | -3 | 10 | 0 | -5 |
| 5 | [12-45] | H | 3 | 1 | 2 | 2 | 1 | 1 |
| | [46-58] | H | 9 | 5 | -2 | 9 | 1 | -5 |
| 6 | [12-49] | H | 3 | 1 | 2 | 2 | 1 | 1 |
| | [50-58] | H | 9 | 5 | -2 | 9 | 1 | -5 |
| 7 | [12-36] | H | 8 | 4 | 1 | 7 | 1 | 0 |
| | [37-58] | H | 9 | 5 | 0 | 9 | 1 | -2 |
| 8 | [12-16] | H | 7 | 3 | 3 | 5 | 2 | 2 |
| | [17-40] | H | 8 | 4 | 1 | 7 | 1 | 0 |
| | [41-58] | H | 9 | 5 | 0 | 9 | 1 | -2 |
| 9 | [12-21] | H | 7 | 3 | 3 | 5 | 2 | 2 |
| | [22-29] | H | 8 | 4 | 1 | 7 | 1 | 0 |
| | [30-58] | H | 12 | 7 | -1 | 12 | 1 | -3 |
| 10 | [12-19] | H | 7 | 3 | 3 | 5 | 2 | 2 |
| | [20-34] | H | 11 | 6 | 1 | 10 | 1 | 0 |
| | [35-58] | H | 12 | 7 | -1 | 12 | 1 | -3 |
| 11 | [12-22] | B | 11 | 4 | 2 | 10 | 1 | 2 |
| | [23] | H | 11 | 5 | 2 | 9 | 2 | 0 |
| | [24-43] | B | 12 | 5 | 1 | 12 | 1 | 0 |
| | [44-58] | H | 12 | 7 | -1 | 12 | 1 | -3 |

**Appendix 2**

Here we give the expanded polynomial form of the formulae for $a_i$ and the cover.

$$a_2 = k_{21}t + c_{21}$$

$$a_3 = k_{32}k_{21}t^2 + (k_{32}c_{21} + c_{32}k_{21} + k_{31})t + c_{32}c_{21} + c_{31}$$

$$\begin{aligned} a_4 = {} & k_{43}k_{32}k_{21}t^3 \\ & + (k_{43}k_{32}c_{21} + k_{43}c_{32}k_{21} + c_{43}k_{32}k_{21} + k_{43}k_{31} + k_{42}k_{21})t^2 \\ & + (k_{43}c_{32}c_{21} + c_{43}k_{32}c_{21} + c_{43}c_{32}k_{21} + k_{43}c_{31} + c_{43}k_{31} + k_{42}c_{21} + c_{42}k_{21} + k_{41})t \\ & + (c_{43}c_{32}c_{21} + c_{43}c_{31} + c_{42}c_{21} + c_{41}) \end{aligned}$$

$$\begin{aligned} C = {} & k_{54}k_{43}k_{32}k_{21}t^4 \\ & + (k_{54}k_{43}k_{32}c_{21} + k_{54}k_{43}c_{32}k_{21} + k_{54}c_{43}k_{32}k_{21} + c_{54}k_{43}k_{32}k_{21} \\ & \qquad + k_{54}k_{43}k_{31} + k_{54}k_{42}k_{21} + k_{53}k_{32}k_{21})t^3 \\ & + (k_{54}k_{43}c_{32}c_{21} + k_{54}c_{43}k_{32}c_{21} + k_{54}c_{43}c_{32}k_{21} + c_{54}k_{43}k_{32}c_{21} + c_{54}k_{43}c_{32}k_{21} + c_{54}c_{43}k_{32}k_{21} \\ & \qquad + k_{54}k_{43}c_{31} + k_{54}c_{43}k_{31} + c_{54}k_{43}k_{31} \\ & \qquad + k_{54}k_{42}c_{21} + k_{54}c_{42}k_{21} + c_{54}k_{42}k_{21} \\ & \qquad + k_{53}k_{32}c_{21} + k_{53}c_{32}k_{21} + c_{53}k_{32}k_{21} \\ & \qquad + k_{54}k_{41} + k_{53}k_{31} + k_{52}k_{21})t^2 \\ & + (k_{54}c_{43}c_{32}c_{21} + c_{54}k_{43}c_{32}c_{21} + c_{54}c_{43}k_{32}c_{21} + c_{54}c_{43}c_{32}k_{21} \\ & \qquad + k_{54}c_{43}c_{31} + c_{54}k_{43}c_{31} + c_{54}c_{43}k_{31} \\ & \qquad + k_{54}c_{42}c_{21} + c_{54}k_{42}c_{21} + c_{54}c_{42}k_{21} \\ & \qquad + k_{53}c_{32}c_{21} + c_{53}k_{32}c_{21} + c_{53}c_{32}k_{21} \\ & \qquad + k_{54}c_{41} + c_{54}k_{41} \\ & \qquad + k_{53}c_{31} + c_{53}k_{31} \\ & \qquad + k_{52}c_{21} + c_{52}k_{21} \\ & \qquad + k_{51})t \\ & + (c_{54}c_{43}c_{32}c_{21} \\ & \qquad + c_{54}c_{43}c_{31} + c_{54}c_{42}c_{21} + c_{53}c_{32}c_{21} \\ & \qquad + c_{54}c_{41} + c_{53}c_{31} + c_{52}c_{21} \\ & \qquad + c_{51}) \end{aligned}$$

**Appendix 3**

This is a transcript of Selmer's notes on "Good Hofmeister Bases" (Selmer, E.S., [8]).

In this transcript I have replaced Greek letters alpha, beta, gamma, delta by A, B, C and D, the summation symbol (Greek capital sigma) by [sigma], the Greek capital "Tor" by T and the Greek capital "Phi" by P. The traditional printer's paragraph sign is replaced by [paragraph], and umlauts are, in general, not mentioned: for example "Braunschadel" should really have umlauts on the "u" and second "a", but neither I nor, it seems, Selmer can be certain or consistent about this.

The first section - pages 1 to 8 - was written in December 1992, probably while Selmer was still at the University; I believe this was before he retired, after which he took a holiday before returning to his home to continue the document (pages 9 to 14) in March 1993 (as if nothing had happened in between!).

[Page 1]

Selmer, Dec. 1992

On "good" Hofmeister bases

The Hofmeister-Schell bases used by Mossige and Challis have the form

$$
(1)\quad \begin{cases} a_2 = (9t + c_{21}) \\ a_3 = (4t + c_{31}) + (3t + c_{32})\, a_2 \\ a_4 = (7t + c_{41}) + (2t + c_{42})\, a_2 + (2t + c_{43})\, a_3 \end{cases}
$$

The corresponding h-range , with

$$h = 12t + r\ ,$$

can be determined by Proposition 6.1 in Mossige's thesis, with my correction of Jan 92. For all bases (1) mentioned below, I have checked the formulas for $n_h(A_4)$ given by Mossige and/or Challis, using Proposition 6.1.

By a "good" basis (1), we mean one with large h-range, often extremal. We have the following good bases available:

(I) Challis' extremal bases (A) in his paper.
(II) Mossige's bases in Table 4 of his thesis, for r = i = 2, 4, 5, 6, 9, 10 (not covered by (I)).
(III) Challis' good bases for 12 <= t <= 20 (his letter to me of April 25, 1992). Mostly covered by (I), and only two essentially new forms:

[Page 2]

| | r | $c_{21}$ | $c_{31}$ | $c_{32}$ | $c_{41}$ | $c_{42}$ | $c_{43}$ | |
|---|---|---|---|---|---|---|---|---|
| (III.1) | 10 | 11 | 6 | 1 | 10 | 1 | 0 | t > 19 |
| (III.2) | 11 | 11 | 5 | 2 | 9 | 2 | 0 | t > 17 . |

The bases at hand have two different forms of $n_h(A_4)$ , types "A" and "B":

$$(2A)\quad n_h(A_4) = (2t + c_{51}) + (t + c_{52})\, a_2 + (6t + c_{53})\, a_3 + (3t + c_{54})\, a_4$$
$$(2B)\quad n_h(A_4) = (3t + C_{51}) + (2t + C_{52})\, a_2 + (4t + C_{53})\, a_3 + (3t + C_{54})\, a_4$$

Both representations have coefficient sum $h = 12t + r$, hence [sigma] c = [sigma] C = r . Type A dominates : All bases (I), bases (II) for r = 4, 5, 6, 9 (values of $c_{5j}$ tabulated by Mossige), and basis (III.2), with

(III,2) $c_{51} = 1,\ c_{52} = 2,\ c_{53} = 1,\ c_{54} = 7$ .

This gives a total of 14 bases - (9 different bases (I)).

In contrast to this, there are only 3 bases
of type B at hand: r = 2, 10 for (II), and (III.1).
We list them all completely:

| | | | $d_1$ | $d_2$ | $d_3$ | $d_4$ | $d_5 - d_6$ | | | | | | |
|---|---|---|---|---|---|---|---|---|---|---|---|---|---|
| | | r | $c_{21}$ | $c_{31}$ | $c_{32}$ | $c_{41}$ | $c_{42}$ | $c_{43}$ | $C_{51}$ | $C_{52}$ | $C_{53}$ | $C_{54}$ | |
| | { (II) | { 2 | 4 | 2 | 0 | 3 | 0 | 0 | -1 | 0 | 0 | 3 | $d_1$-2 |
| (3) | { (II) | { 10 | 10 | 5 | 2 | 8 | 1 | 2 | 1 | 1 | 4 | 4 | |
| | { (III.1) | 10 | 11 | 6 | 1 | 10 | 1 | 0 | 1 | 1 | 0 | 8 | |

The $c_{5j}$ of (III.2) and $C_{5j}$ of (III.1) (not given
by Challis) were calculated by me from Proposition 6.1.

[Page 3]

So far, I have only systematized
the "good" Hofmeister bases we have at
hand. Now, I shall mention some very
striking observations I have made for
such bases. They resulted when I trans-
formed the representations (2A), (2B) to regu-
lar form .

For (2A), the corresponding regular repre-
sentation in all 14 cases is

$$\begin{aligned}(4A)\quad n_h(A_4) = &\ 5t + c_{51} + (c_{21} + 2c_{31} - 2c_{41}) \\ &+ \{ 3t + c_{52} + (2c_{32} - 2c_{42} - 1) \} a_2 \\ &+ \{ 2t + c_{53}) - (2c_{43} + 2) \} a_3 \\ &+ \{ 3t + c_{54} + 2) \} a_4\end{aligned}$$

The coefficient sum is 13t + $r_1$, where in
fact $r_1$ = r in most cases (but 4 exceptions).
We note that when using Proposition 6.1
- which always gives the regular form -
we must then end up with $n_2(A_4)$ in (6.5).

Denoting $n_h(A_4)$ of (4A) by

$$n_h(A_4) = (5t + d_{51}) + (3t + d_{52}) a_2 + (2t + d_{53}) a_3 + (3t + d_{54}) a_4 ,$$

there are two conditions for regularity:

$$(5t + d_{51}) + (3t + d_{52}) a_2 < a_3$$
$$(5t + d_{51}) + (3t + d_{52}) a_2 + (2t + d_{53}) a_3 < a_4 .$$

With the right hand sides given by (1), we

[Page 4]

get the two necessary and sufficient conditions (at least for sufficiently large t):

$d_{52} < c_{32}$ , $d_{53} < c_{43}$ ,

which take the form

$c_{52} <= 2c_{42} - c_{32}$ , $c_{53} <= 3c_{43} + 1$ .

When checking whether this was satisfied for the bases of type A, I discovered that in all 14 cases, these conditions were satisfied with equality! It was then natural to check if a similar linear relation exists between the $c_{21}$, and in fact $3c_{31} - 2c_{41} + c_{51} + 2 = 0$ holds in all cases.

It thus seems that for our bases (1) of type A, the h-range (2A) can be immediately determined by

(5A) $c_{51} = 2c_{41} - 3c_{31} - 2$ , $c_{52} = 2c_{42} - c_{32}$ , $c_{53} = 3c_{43} + 1$, [ this line is strongly emphasised with ] [ multiple vertical lines in the left hand ] [ margin ]

together with the obvious relation [sigma] $c_{5j} = r$ :

$c_{54} = r - (c_{51} + c_{52} + c_{53})$ .

Substituting (5A) in (4A) , we get

(6A) $n_h(A_4) = (5t + c_{21} - c_{31} - 2) + (3t + c_{32} - 1) a_2 + (2t + c_{43} - 1) a_3 + (3t + c_{54} + 2) a_4$ .

This form is then equivalent to the observation (5A) .

[Page 5]

Now we can prove (6A) and thus (5A), using Mossige's Proposition 6.1. His formulas (6.5) give the regular form of $n_h(A_4)$ in the different cases, and should therefore be compared with (4A) above. As already mentioned, Mossige's form $n_2(X,P)$ is the only one giving the correct t-coefficients of (4A) (t = j in his notation). If we substitute in $n_2(X,P)$ the coefficients

$d_1 = c_{21}$ ; $d_2 = c_{31}$ , $d_3 = c_{32}$ ; $d_4 = c_{41}$ , $d_5 = c_{42}$ , $d_6 = c_{43}$ ,

we just get the three first terms in (6A)

(the coefficients of $a_4$ cannot be compared).
This completes the promised proof.

If I had discovered directly
the above proof of (6A), I would in all
probability not have performed the calcu-
lation back to the interessting observa-
tion (5A).

One final remark: In my note of Jan. 92,
the corrections to Mossige's Proposition 6.1
did not affect the function $f_2(X,P)$ of (6.3),
corresponding to $n_2(X,P)$. But my condition for
a necessary correction:

(7) $d_1 + 3d_2 - 3d_4 >= 0$ or $2d_3 - 3d_5 - 2 >= 0$ ,

are satisfied in some of the 14 cases of type A,
and not in other cases.

[Page 6]

We now turn to the "good" Hofmeister
bases of type B, with h-range given by
(2B). We only have the three bases (3) at hand.

The corresponding regular representation in
these cases is

(4B)
$$\begin{aligned} n_h(A_4) = \; & 9t + C_{51} + (c_{21} + c_{31} - c_{41}) \\ & + \{3t + C_{52} + (c_{32} - c_{42} - 1)\}\, a_2 \\ & + \{2t + C_{53} - (c_{43} + 1)\}\, a_3 \\ & + \{3t + C_{54} + 1\}\, a_4 \end{aligned}$$

$$= 9t + D_{51} + (3t + D_{52})\, a_2 + (2t + D_{53})\, a_3 + (3t + D_{54})\, a_4 \quad \text{(say)}.$$

The necessary and sufficient conditions for
regularity are now

$$D_{52} < c_{32} \, , \; D_{53} < c_{43} \, ,$$

which take the form

$$C_{52} <= c_{42} \, , \; C_{53} <= 2c_{43} \, .$$

Again, these are satisfied with equality for
the three bases (3):

(5B) $C_{52} = c_{42}$ , $C_{53} = 2c_{43}$ .

But to find an expression for $C_{51}$ , in analogy
with $c_{51}$ of (5A), is difficult from the very few bases (3) at hand.

Substituting (5B) in (4B), we get

(6B) $n_h(A_4) = (9t + C_{51} + c_{21} + c_{31} - c_{41}) + (3t + c_{32} - 1)\, a_2$
$+ (2t + c_{43} - 1)\, a_3 + (3t + C_{54} + 1)\, a_4$ .

[Page 7]

For the corresponding (6A), we could identify this uniquely with the function $n_2(X,P)$ of (6.5) in Mossige , by looking at the t-coefficients. For (6B), however, these coefficients leave both $n_3(X,P)$ and $n_4(X,P)$ as candidates. For the bases (3), only $n_3(X,P)$ (possibly modified, cf. below) turns up. This was to be expected, since both coefficients of $a_2$ and $a_3$ in (6B) fit for $n_3(X,P)$ but not for $n_4(X,P)$.

But now another complication turns up, because of my correction of Mossige's Proposition 6.1. If the condition (7) is not satisfied, Prop. 6.1 goes unmodified, with the constant term

$$9t + 2d_1 + 3d_2 - 3d_4 - 2$$

of $n_3(X,P)$. A comparison with (6B) shows that we must then expect

$$C_{51} = c_{21} + 2c_{31} - 2c_{41} - 2 .$$

This holds for the last basis (3), which is the only one where (7) fails. The numerical material is not overwhelming!

If on the other hand (7) holds, we must correct Proposition 6.1 as in my note of Jan. 92. In particular, $n_3(X,P)$ is replaced by $n_{3B}(X,P)$, with a different constant term

$$9t + d_1 - 2 .$$

[Page 8]

We must then expect

$$C_{51} = c_{41} - c_{31} - 2 ,$$

and this holds for the first two bases (3) (which both satisfy (7)).

Because of these complications, it may
safely be said that the (dominating)
type A is the most interesting one!

I have asked earlier why the
good Hofmeister bases seem to give
larger h-ranges than the Braunshadel
bases. Now a similar question turns
up, why type A is usually "better"
than type B.

It seems very difficult to answer any
of the questions. But one thing would be
interesting to look into: In his good
bases (III), Challis list one (III.1) of type B as
high up as for t > 19. I therefore have
the following question to him: Would it take [ this and the following two lines are strongly ] [ emphasised with multiple vertical lines ] [ in the left hand margin ]
much time to extend the list of "good"
Hofmeister bases to t > 20?

So far, all examined (by Prop. 6.1) Hofmeister [ this entire paragraph is emphasised as ] [ as described above ]
bases have led to $n_2(X,P)$ or $n_3(X,P)/n_{3B}(X,P)$.
It would be interesting to give some bases
where other functions $n_l(X,P)$ turn up, that is,
bases which are not of type A or B. Preferably,
such bases should be "good".

[Page 9]

Selmer, March 1993

On "good" Hofmeister bases
(continued).

I have now received Challis' note
"Looking for good Hofmeister and Braun-
schadel bases" (17 pp, but not paginated!),
with very many interesting observations.

I asked about the occurrence of my
"dominating" type A, with $k_{51} - k_{54} = (5, 3, 2, 3)$
by (4A), versus type B with (9, 3, 2, 3) in (4B).
In Challis' Table 4, my basis (III.1) is
- correctly - the only case B (for r = 10). More
interesting is of course Challis' extension
to t <= 58 in Table 10 of "best" Hofmeister
bases. It still contains only A and B, but
now with four cases B. However, only
two of these, for r = 0 and r = 4, apply

for the largest t = 58, cf. Table 8. So
it seems that type A is "strongly" dominating for "very good" Hofmeister bases for large t. For the "not so good" bases for r = 0 in Table 7, only half (14 of 29) of the bases are of type A. Quite surprisingly, type B has disappeared completely, but there are three other types (6, 3, 2, 3), (8, 1, 2, 3) and (5, 2, 2, 3) instead. The first of these inspired Challis' correction of Mossige's Proposition 6.1.

[Page 10]

For Braunschadel bases, a similar dominance appears for $k_{51}$ - $k_{54}$ = (4, 3, 2, 3), at least for larger t in Table 11. This now contains seven bases of other types, but again only two of these, for r = 0 and r = 11, apply for t = 58, cf. Table 9.

In [para] 5, "The future", Challis conjectures that all extremal bases for sufficiently large h are either Hofmeister or Braunschadel - and the coefficient sets $c_{ij}$, for a fixed value of r, will "move slowly" away from (0, 0, 0, 0, 0, 0) without limit.

Again, I must return to the question: What is really (and "naturally") a Hofmeister or Braunschadel basis? In my "Comments to Challis' letter of 19.10.92" (Dec. 1992), I exemplify this question with the bases (1) and (2) page C, of which (1) is not Hofm/Braunsch. A similar and much more striking example will appear below.

Challis also wonders whether one or the other "form" will eventually dominate. I assume that by form, he refers to Hofmeister versus Braunschadel. But we may also ask if, among for instance extremal Hofmeister bases, there will be a dominance of "types" (like A above), that is, of the choice $n_i$ in Mossige's (6.5),

[Page 11]

with Challis' latest addition.

For the determination of $n_i$, we cannot use only (asymptotic) parameter bases, but need exact coefficients.

Outside Challis' range t <= 58, we only have one such case at disposal, namely Mossige's "optimal" basis of Theorem 11.1 (in thesis, = Th. 6.1 in Math. Scand paper). Putting bt = T, and then replacing At by t, this basis may be written as

$$
(8)\quad \begin{cases} h = 12t \quad (r = 0) \\ a_2 = 9t + 15T \\ a_3 = 4t + 14T + (3t - 15T + 2)\, a_2 \\ a_4 = 9t + 23T + (2t - 2T)\, a_2 + (2t - 20T)\, a_3 . \end{cases}
$$

From Mossige's condition A >= 25b, we get T <= t/25.
Is this a Hofmeister basis or not? [ this line is emphasised as described previously]
Judging from the t-coefficients (9, 4, 3, 7, 2, 2) , one would say "yes". On the other hand, Mossiges very good choice T = t/206 (see also his example p. 48) gives a non-Hofmeister basis with prefactor P = 2.008.

But under all circumstances, the h-range of the basis (8) can be calculated from Mossige's Proposition 6.1, which was specially designed to find the h-range of a Hofmeister basis. In (6.3), we must use $f_{1A}$ , and $2d_3 - 2d_5 - 2 < 0$, so Challis'

[Page 12]

latest modification does not apply, nor does my earlier modification. In Prop. 6.1 , we then find Z = 45T - 1, L = {1, 7}, and the h-range min $\{n_1, n_7\} = n_1$ is just the one given by Theorem 11.1.

The "type" $k_{51} - k_{54} = (8, 1, 0, 3)$ is not the type A which dominates for smaller t. However we shall see that it is easy to construct alternative "optimal" bases of several different types, more specifically corresponding to $n_2$ (type A), $n_3$ , $n_6$ and $n_7$.

Using Prop. 6.1 on the basis (8), we find

$$
(9)\quad \left.\begin{matrix} f_{1A} \\ f_2 \\ f_3 \\ f_6 \\ f_7 \end{matrix}\right\} = -X + 45T + \begin{cases} -1 \\ 4 \\ 2 \\ 2 \\ -1 \end{cases} \qquad \begin{matrix} f_4 = -X + 57T \\ f_5 = -X + 51T - 1 \end{matrix}
$$

The constant addends to the left result from
the choice P = (0, 0, 2, 0, 0, 0) in Theorem 11.1.
The smallest addends -1 determine L = {1, 7}
in line 3 above. By changing P, we can
get other mininal addends, for instance

P = ( 1, 0, 0, 0, 0, 0) : min for $f_3$ }
P = (-1, 0, 0, 0, 0, 0) : min for $f_6$ }   Z = 45T + 1 in all cases
P = (0, 0, 0, 0, -1, 0) : min for $f_7$ }

[Page 13]  [ and Page 14 ]

So the types corresponding to $n_3$, $n_6$
and $n_7$ all turn up.

But the most interesting - and
"good looking" - choice is <u>all zeros</u> ,
P = (0, 0, 0, 0, 0, 0). Then Z = 45T + 2, L = {2, 3, 6} ,
and the minimum at the bottom of Mossige's
p. 39 occurs for $n_2$ (type A), with

(10)  {   $n_h(A_4) = (3t + 45T + 3)\ a_4 + (2t - 20T -1)\ a_3$
{                        $+ (3t - 15T + 1)\ a_2 + (5t + T - 2)$ .

Of course, all such cases contain
the same term (3t + 45T) $a_4$, and thus
give the <u>same asymptotic value</u> of $n_h(A_4)$.
The advantage of the last, "all zero" case
is that it is <u>analagous to my asymptotic
parameter bases</u>, where I drop all constant
terms. For instance, the addend +2 in
the $a_2$-coefficient in Mossige's example p. 48
now disappears, resulting in the h-range

$$n(h, C) = (663t + 3)\ a_4 + (392t - 1)\ a_3 + (603t + 1)\ a_2 + (1031t - 2)$$

of (10).

When I asked Mossige why he chose
the vector P = (0, 0, 2, 0, 0, 0), he explained
that it was to get an <u>$h_0$-basis</u>, since
his condition of (6.4),

$$f_8(P) = - d_1 - d_3 + i + 2 >= 0 \quad (i = 0)$$

is then <u>sharp</u>. With the zero vector,
$f_8(P) = 2$, meaning that <u>$h = h_0 + 2$</u>.

It follows from (9) that with the
same vector B, the $A_4$ of Theorem 11.1
can never have a h-range of the form

("type") $n_4$ or $n_5$, whatever the choice of P. This raises an interesting (and certainly difficult) question: Can such a h-range, with the same asymptotic prefactor 2.008, be obtained from a different choice of the vector B ?
Generally, is Mossige's B "God-given"?   [ this line is emphasised as described previously ]

[The End]

**References**

[1]  Braundshadel, R., "Zum Reichweitenproblem", Staatsexamensarbeit, Math. Inst., Joh. Gutenberg-Univ., Mainz (1987).

[2]  Braundshadel, R., "Zum Reichweitenproblem", Diplomarbeit, Math. Inst., Joh. Gutenberg-Univ., Mainz (1988).

[3]  Challis, M.F., "Two new techniques for computing extremal h-bases $A_k$", Computer J. 36 (1993), pp117-126.

[4]  Challis,M.F., "The Postage Stamp Problem: Formulae and proof for the case of three denominations", Storey's Cottage, Whittlesford, Cambridge (1990).

[5]  Challis, Michael F. & Robinson, John P., "Some Extremal Postage Stamp Bases", Journal of Integer Sequences, Vol 13 (2010), Article 10.2.3.

[6]  Mossige, Svein, "On the Extremal h-range of the Postage Stamp Problem with Four Stamp Denominations", Research Monograph no. 41-06-01-86, ISSN 0332-5407, Department of Mathematics, University of Bergen, Norway (1986).

[7]  Mossige, Svein, "On extremal h-bases $A_4$", Math Scand. 61 (1987).

[8]  Selmer, E. S., "On good Hofmeister Bases", Private Communication (1992/1993).

[9]  Selmer, Ernst S., "Asymptotic h-ranges and Dual Bases", Research Monograph no. 56-02-02-90, ISSN 0332-5407, Department of Mathematics, University of Bergen, Norway (1990).

[10]  Selmer, Ernst S, and Mossige, Svein, "Stohr sequences in the postage stamp problem", Research Monograph no. 32, ISSN 0332-5407, Department of Mathematics, University of Bergen, Norway (1984).

[11]  Selmer, E. S., "The Local Postage Stamp Problem Part 1: General Theory", Research Monograph no. 42-04-15-86, ISSN 0332-5407, Department of Mathematics, University of Bergen, Norway (1986).

[12]  Selmer, E. S., "The Local Postage Stamp Problem Part 2: The Bases $A_3$ and $A_4$", Research Monograph no. 44-09-15-86, ISSN 0332-5407, Department of Mathematics, University of Bergen, Norway (1986).

[13]  Selmer, E. S., "The Local Postage Stamp Problem Part 3: Supplementary Volume", Research

Monograph no. 57-06-12-90, ISSN 0332-5407, Department of Mathematics, University of Bergen, Norway (1990).

[14] Selmer, E.S., "On the Postage Stamp Problem with three stamp denominations", Math Scand. 47 (1980), pp 29-71.

[15] Selmer, E.S. and Rodne, A., "On the Postage Stamp Problem with three stamp denominations II", Math Scand. 53 (1983), pp 145-156.

[16] Selmer, E.S., "On the Postage Stamp Problem with three stamp denominations III", Math Scand. 56 (1985), pp 105-116.